\documentclass[12pt]{amsart}
\usepackage{amsthm}
\usepackage{amssymb}
\usepackage{amsbsy}
\usepackage{amscd}
\usepackage[mathscr]{eucal}
\usepackage{nicefrac}

\DeclareSymbolFont{rsfs}{U}{rsfs}{m}{n}
\DeclareSymbolFontAlphabet{\mathrsfs}{rsfs}
\usepackage{a4wide}
\usepackage[all]{xy}
\usepackage{tikz}
\usetikzlibrary{arrows,decorations.pathmorphing,backgrounds,positioning,fit,petri}
\usetikzlibrary{decorations.markings}

\usepackage{verbatim}
\usepackage{version}
\usepackage{color}
\newenvironment{NB}{
\color{red}{\bf NB}. \footnotesize
}{}
\excludeversion{NB}
\excludeversion{NB2}
\makeatletter

\usepackage{url}
\usepackage{ifmtarg}

\usepackage{xr-hyper}
\usepackage[
CJKbookmarks,
pdftex,
bookmarks=true,
colorlinks=true,
pdfnewwindow=true]{hyperref}
\usepackage{cleveref}
\let\RR\relax
\let\bN\relax

\let\cO\relax
[affine_pre.pdf]
[http://arxiv.org/pdf/1604.03625.pdf]
\usepackage{xspace}

\usepackage[all]{xy}

\usepackage{stmaryrd}
\SetSymbolFont{stmry}{bold}{U}{stmry}{m}{n}
\usepackage{pigpen} % for pigpenfont used for fiber squares
\usepackage{enumitem}
\usepackage{calc}

\usepackage{textgreek}

\usepackage{CJK}
\renewcommand{\thesubsection}{\thesection(\@roman\c@subsection)}
\newcounter{number}
\makeatother
\begin{CJK}{UTF8}{ipamp}
\newtheorem{Theorem}[equation]{定理}
\end{CJK}

\theoremstyle{definition}

\theoremstyle{remark}

\numberwithin{equation}{section}
\usepackage{cleveref}
\crefname{Theorem}{Theorem\xspace}{Theorems}
\crefname{section}{\S}{\S\S}
\crefformat{section}{\S#2#1#3} % to remove unnecessary space
\crefmultiformat{section}{\S\S#2#1#3}{, #2#1#3}{, #2#1#3}{, #2#1#3}
\crefname{Lemma}{Lemma\xspace}{Lemmas\xspace}
\crefformat{Lemma}{Lemma~#2#1#3}
\crefname{Proposition}{Proposition\xspace}{Propositions\xspace}
\crefformat{Proposition}{Proposition~#2#1#3}
\crefname{Corollary}{Corollary\xspace}{Corollaries\xspace}
\crefformat{Corollary}{Corollary~#2#1#3}
\crefname{Definition}{Definition}{Definitions}
\crefformat{Definition}{Definition~#2#1#3}
\crefname{Remark}{Remark\xspace}{Remarks\xspace}
\crefformat{Remark}{Remark~#2#1#3}
\crefname{Remarks}{Remark\xspace}{Remarks\xspace}
\crefformat{Remarks}{Remark~#2#1#3}
\crefname{Conjecture}{Conjecture\xspace}{Conjectures\xspace}
\crefformat{Conjecture}{Conjecture~#2#1#3}
\crefname{figure}{Figure\xspace}{Figure\xspace}
\crefformat{figure}{Figure~#2#1#3}

\crefformat{appendix}{\S#2#1#3}
\crefformat{enumi}{#2#1#3}

\newcommand{\secref}[1]{\S\ref{#1}}
\newcommand{\subsecref}[1]{\S\ref{#1}}
\newcommand{\CC}{{\mathbb C}}
\newcommand{\ZZ}{{\mathbb Z}}

\newcommand{\RR}{{\mathbb R}}
\newcommand{\SU}{\operatorname{\rm SU}}
\newcommand{\GL}{\operatorname{GL}}

\newcommand{\SO}{\operatorname{\rm SO}}
\newcommand{\grpSp}{\operatorname{\rm Sp}}

\newcommand{\su}{\operatorname{\mathfrak{su}}}

\newcommand{\Spec}{\operatorname{Spec}\nolimits}

\newcommand{\Hom}{\operatorname{Hom}}

\renewcommand{\MR}[1]{}

\newcommand{\dslash}{/\!\!/}
\newcommand{\vin}[1]{\operatorname{i}(#1)} % incoming vertex
\newcommand{\vout}[1]{\operatorname{o}(#1)} % outgoing vertex
\newcommand{\bM}{\mathbf M}
\newcommand{\bN}{\mathbf N}

\newcommand{\tslash}{/\!\!/\!\!/}
\newcommand{\tslabar}{\mathbin{
\setbox0=\hbox{/\!\!/\!\!/}\rule[0.4\ht0]{\wd0}{.3\dp0}\kern-\wd0\box0}}

\newcommand{\la}{\lambda}

\newcommand{\Gr}{\mathrm{Gr}}

\newcommand{\cR}{\mathcal R}
\newcommand{\cF}{\mathcal F}
\newcommand{\cT}{\mathcal T}
\newcommand{\cK}{\mathcal K}
\newcommand{\cO}{\mathcal O}

\newcommand{\scP}{\mathscr P}

\newcommand{\St}{\mathrm{St}}

\newcommand{\Uh}[2][G]{{\mathcal U}_{#1}^{#2}}
\newcommand{\Bun}[2][G]{\operatorname{Bun}_{#1}^{#2}}

\makeatletter
\newcommand{\cA}[1][{}]{%
  \@ifmtarg{#1}%
  {\mathcal A}% if #1 is empty
  {\mathcal A(#1)}% if #1 is not empty
}
\newcommand{\cAh}[1][{}]{%
  \@ifmtarg{#1}%
  {\mathcal A_\hbar}% if #1 is empty
  {\mathcal A_\hbar(#1)}% if #1 is not empty
}
\makeatother

\newcommand{\DD}{\mathbb D}
\newcommand{\DC}{\boldsymbol\omega}
\newcommand{\Stab}{\operatorname{Stab}}

\newcommand{\ft}{\mathfrak t}

\newcommand{\po}{\ar@{}[dr]|{\text{\pigpenfont R}}}
\newcommand{\pb}{\ar@{}[dr]|{\text{\pigpenfont J}}}
\newcommand{\pp}{\ar@{}[dr]|{\text{\pigpenfont P}}}

\newcommand{\cM}{\mathcal M}

\makeatletter
\AtBeginDvi{\pdfmapline{=ipamp@Unicode@  <ipamp-select.ttf}}
\DeclareFontFamily{C70}{ipamp}{\hyphenchar \font\m@ne}
\DeclareFontShape{C70}{ipamp}{l}{n}{ <-> CJK * ipamp}{}
\DeclareFontShape{C70}{ipamp}{m}{n}{ <-> CJK * ipamp}{\CJKnormal}
\DeclareFontShape{C70}{ipamp}{bx}{n}{ <-> CJKb * ipamp}{\CJKbold}
\makeatother

\begin{document}

\title[Coulomb branches of $3d$ $N=4$ gauge theories] {Introduction to a
provisional mathematical definition of Coulomb branches of $3$-dimensional
  $\mathcal N=4$ gauge theories
%{\rm Preliminary Version (\today)}
}
% \author[A.~Braverman]{Alexander Braverman}
% \address{
% Department of Mathematics,
% University of Toronto and Perimeter Institute of Theoretical Physics,
% Waterloo, Ontario, Canada, N2L 2Y5
% }
% \email{braval@math.toronto.edu}
% \author[M.~Finkelberg]{Michael Finkelberg}
% \address{National Research University Higher School of Economics,
% Department of Mathematics, 20 Myasnitskaya st, Moscow 101000 Russia}
% \email{fnklberg@gmail.com}
\author[H.~Nakajima]{Hiraku Nakajima}
%\begin{CJK}{UTF8}{ipaexm}
%\author[
% %\begin{CJK}{UTF8}{ipaexm}
% \authorname %中島　啓
% %\end{CJK}
% ]
% {%\begin{CJK}{UTF8}{ipaexm}
% \authorname %中島　啓
% %\end{CJK}
%  (Hiraku Nakajima)}
%\author[中島　啓]{中島　啓}
% \address{Research Institute for Mathematical Sciences,
% Kyoto University, Kyoto 606-8502,
% Japan}
\address{\begin{CJK}{UTF8}{ipamp}
〒606-8502　京都市左京区北白川追分町
　京都大学
　数理解析研究所
\end{CJK}}
\email{nakajima@kurims.kyoto-u.ac.jp}

\subjclass[2000]{}

\renewcommand{\abstractname}{\begin{CJK}{UTF8}{ipamp}
概要\end{CJK}
}

\begin{abstract}
\begin{CJK}{UTF8}{ipamp}
  この論説は、\cite{2015arXiv150303676N,main}で提唱された、$3$次元$\mathcal N=4$超対称性ゲージ理論のクーロン枝の、暫定的な数学的定義に関する入門である。
\end{CJK}

  This is an introduction to a provisional mathematical definition of
  Coulomb branches of $3$-dimensional $\mathcal N=4$ supersymmetric
  gauge theories, studied in \cite{2015arXiv150303676N,main}.
\end{abstract}

\maketitle

\setcounter{tocdepth}{2}
%\tableofcontents

\begin{NB}
\section*{Some macros}
\begin{verbatim}
\newcommand{\Gr}{\mathrm{Gr}}
\newcommand{\cR}{\mathcal R}
\newcommand{\cF}{\mathcal F}
\newcommand{\cT}{\mathcal T}
\newcommand{\cK}{\mathcal K}
\newcommand{\cO}{\mathcal O}

\newcommand{\scP}{\mathscr P}


\newcommand{\St}{\mathrm{St}}

\newcommand{\Uh}[2][G]{{\mathcal U}_{#1}^{#2}}
\newcommand{\Bun}[2][G]{\operatorname{Bun}_{#1}^{#2}}
\newcommand{\DD}{\mathbb D}
\newcommand{\DC}{\boldsymbol\omega}
\newcommand{\Stab}{\operatorname{Stab}}
\end{verbatim}

\verb+\cAh+ yields $\cAh$.

\verb+\cAh[G,\bN]+ yields $\cAh[G,\bN]$.

\verb+\cAh[G]+ yields $\cAh[G]$.

\verb+\cAh[]+ yields $\cAh[]$.

\end{NB}

\begin{CJK}{UTF8}{ipamp}
%\begin{CJK}{UTF8}{min}

%\input{introduction}

\section{複素シンプレクティック多様体と変形量子化}

$G$を複素簡約群とし、$\bM$をそのシンプレクティックな表現とする。すなわち、$\bM$は$\CC$上のシンプレクティック形式を持つベクトル空間であり、$G$はシンプレクティック形式を保って線形に作用している。

%物理的な背景を忘れて、現時点で数学的に厳密に確立されている部分だけを見ると、
%
{\bf Coulomb枝}$\cM_C\equiv\cM_C(G,\bM)$は、理論物理における場の量子論の研究に動機付けられて発見された、$(G,\bM)$からアファイン複素シンプレクティック多様体\footnote{一般には特異点を持つ。高々Beauvillの意味でシンプレクティックな特異点しか持たないと期待されているが、証明は与えられていない。}を作るレシピである。
\begin{equation*}
  (G,\bM) \leadsto \cM_C(G,\bM)
\end{equation*}
作り方は、これまで知られている代数多様体の与え方、多項式の零点、商空間、等々とはかなり毛色が異なる。まず座標環$\CC[\cM_C]$を、幾何学的表現論でよく使われるホモロジー群とその上の合成積を考える方法で作る。そしてその可換環のスペクトラムとして$\cM_C$を定め、その幾何学的な性質を調べる、という手法を取る。

あとで説明するように、$\cM_C$は $T^* T^\vee/W$と双有理同型である。
\begin{equation*}
    \cM_C \approx T^*T^\vee/W = \ft\times T^\vee/W
\end{equation*}
ここで、$T^\vee$は$G$の極大トーラス$T$の双対トーラスであり、$W$はワイル群である。$T^*T^\vee$は\linebreak[3]$T^\vee$のco\-tangent bundleで、$\ft$は$T$のリー環である。
特に、$\cM_C$の双有理類は表現$\bM$には依存しない。

上で言及したように、ホモロジー群とその上の合成積を用いて環を作るレシピは、幾何学的表現論ではよく使われてきた。表現論の研究が目的であるから、そこでは、非可換環を構成するのが普通である。
実際、Coulomb枝においてもその構成法から、$\cM_C$の変形量子化$\cAh$が同時に作られる。ここで、変形量子化とは、$\CC[\hbar]$上で定義された非可換環 $\cAh$ であって、$\cAh/\hbar\cAh$が
$\cM_C$の座標環 $\CC[\cM_C]$ に等しく、Poisson 括弧
\begin{equation*}
    \{ f, g\} = \left.
    \frac{\tilde f\tilde g - \tilde g\tilde f}{\hbar}\right|_{\hbar=0}, \qquad
    \tilde f|_{\hbar=0} = f, \quad \tilde g|_{\hbar=0} = g
\end{equation*}
が、シンプレクティック構造から来るものに一致しているものをいう。これを{\bf 量子化されたCoulomb枝}とよぶ。

振り返って考えれば、表現論で研究されてきた非可換環は、可換環の変形として得られているものが多い。しかし、合成積を用いて可換環を新しく系統的に構成しようという発想は、今回の研究で初めて現れたものだと認識している。

最初の論文\cite{2015arXiv150303676N}では、一般の$\bM$を考えていたが$\CC[\cM_C]$のベクトル空間としての構成にとどまり、積の定義はあとの\cite{main}で与えられた。その際に、$\bM = \bN\oplus\bN^*$という形であると仮定した。この仮定は技術的なものなのか、もしくはより本質的なものなのかはまだ分からないが、物理ではこの仮定を満たさないCoulomb枝も考察されており、どんな条件が満たされていれば定義ができるのか、検討の余地が残されている。なお$\bM=\bN\oplus\bN^*$を仮定する次々節以降では、$\cM_C(G,\bN)$という記号を使うが、混乱のおそれはないと思われる。

$(G,\bM)$ に対して、アファイン複素シンプレクティック多様体を与える、よく知られたレシピがある。それは、シンプレクティック商
\begin{equation*}
    \bM\tslash G = \mu^{-1}(0)\dslash G
\end{equation*}
であり、物理では{\bf Higgs枝}とよばれる。Coulomb枝に対応して$\cM_H \equiv \cM_H(G,\bM)$であらわす。
上の式で、$\mu\colon\bM\to\operatorname{Lie}G^*$は運動量写像であり、$\mu^{-1}(0)$を$G$で(幾何学的不変式論の意味で)割ってできる商空間が$\mu^{-1}(0)\dslash G$である。

$\bM=\bN\oplus\bN^*$となっているときには、$\bN$の上の多項式係数の微分作用素の全体のなす非可換環$\mathcal D(\bN)$が、$\bM$の変形量子化になる。($\hbar$を入れるには、次数によるフィルターに関してRees代数を作る。) シンプレクティック商と同様に、$\mathcal D(\bN)$の$G$作用に関する`商'を作る構成法が、量子シンプレクティック簡約として知られており、それが $\cM_H$の変形量子化を与える。

Higgs枝として、箙多様体やトーリック超ケーラー多様体を例として、表現論的に興味深いシンプレクティック多様体や、その変形量子化が現れることを経験している。
一方、Coulomb枝の研究は始まったばかりであるが、Higgs枝としては得られないシンプレクティック多様体(正確には、有限次元のシンプレクティック・ベクトル空間のシンプレクティック商としての記述が知られていない空間)も現れるので、今後重要性が高まるレシピであると期待している。

また、同じ $(G,\bM)$ からできるHiggs枝とCoulomb枝は、Braden-Licata-Proudfoot-Webster \cite{2014arXiv1407.0964B}の意味で、シンプレクティック双対であることが期待されている。シンプレクティック双対は、複素シンプレクティック多様体のペアの間に不可思議な関係があることを期待するもので、全体像はまだまだ研究途中で見えていないが、少なくともHiggs枝とCoulomb枝の両方を同時に研究することに意味があり、重要であることを示唆している。
\cite{2014arXiv1407.0964B}においては、そのような複素シンプレクティック多様体のペアの例が例示されていたにとどまっていたが、Coulomb枝による系統的な構成が与えられたことになる。ただし、\cite{2014arXiv1407.0964B}で期待されている不可思議な関係のチェックは、今後の課題である。特に、\cite{2014arXiv1407.0964B}は定式化において二つの複素シンプレクティック多様体は、ともにシンプレクティックな特異点解消を持つことが仮定されていたが、多くのHiggs枝、Coulomb枝においてこの仮定は成立しないので、何を期待するのか、ということまで含めて検討する必要がある。

\section{物理的な背景}

前節の説明で、Coulomb枝の数学的な研究に意味があることが伝えられたと期待するが、今節では物理的な背景について、筆者の理解できる範囲内で説明を試みる。ここに書いてあることを理解する必要はないし、筆者自身もよく理解したとは思っていないが、次節で説明する定義がどこから発見されたのかを理解するためと、今後新たな研究成果を上げるためには、背景にある物理のある程度の理解が必要であろうと思っている。

先を急ぐ読者は、この節を飛ばして読んでも構わないが、より深い理解を求める方は、今節を読み、また物理の文献に挑戦していただきたい。物理の文献は\cite{2015arXiv150303676N}にあげたので、これを参照すること。

また、この論説は\cite{main}と同様に、この節以外は物理を知らなくても読めるように書かれており、物理の文献の引用をしない。これは、あくまで読みやすさのためのもので、原典は\cite{2015arXiv150303676N}の文献表から見つけてあたってもらいたい。

物理では、微分幾何と同様に複素簡約群 $G$ の代わりに、その極大コンパクト部分群$G_c$ を取り扱う。同様に、$\bM$には$G_c$で保たれる内積が入っているものとする。

組$(G_c,\bM)$に対して、物理学者は$3$次元の$\mathcal N=4$超対称ゲージ理論を定める。
これは、場の全体のなす無限次元の空間の上の汎関数(ラグランジアン)を与え、量子化して得られる場の量子論の例である。
場のうちで主要なものは、$\RR^3$上の$G_c$主束$P$上の接続と、$\bM$に値を持つような$P$の切断である。さらにいくつかのベクトル束の切断を場に加える。足される場は物理的には必要であるが、ここでは雑な理解しか与えないので、説明は省く。いずれにせよ、物理学者は、接続の曲率や切断の微分を含んだラグランジアンを書き、場の量子論を考える。
ラグランジアンが極小値を取るような接続や切断(と説明を省略した場)のconfigurationは、量子力学でいうところの古典解に対応しており、基本的な対象である。
今の状況では、極小値をとる場は、ただ一つではなく、有限自由度を持った空間になっている。これは、物理では{\bf 真空のモジュライ空間}とよばれる。
%
%また、超対称性を持った場の理論であることから、無限次元の場の全体のなす空間上のFeynmann積分が、真空のモジュライ空間上の有限次元の積分に等しいことが、期待される。

上に述べたように、ラグランジアンは接続の曲率や切断の微分などの和として与えられる。極小値を与える場は、和のうちのいくつかの項が消えているものであり、どの項が消えているかで分けて、真空の枝という考え方をする。その中の典型的なものがHiggs枝$\cM_H$と\linebreak[3]Coulomb枝$\cM_C$である。
Higgs枝$\cM_H$は前節に述べたシンプレクティック商であり、微分幾何学的には超ケーラー商である。接続は自明接続で、定数な切断だけが生き残るので、$\bM$の情報だけが残って超ケーラー商になる。ここでは、超ケーラー商の定義は復習しないので、例えば \cite{MR1193019}を参照してほしい。
例えば筆者が長年に渡り研究している箙多様体や、トーリック超ケーラー多様体は、Higgs枝の例になっている。

一方、Coulomb枝は $(T^\vee_c\times(\RR^3\otimes\ft_c))/W$ となる。$T^\vee_c$は、$G_c$の極大トーラス $T_c$ の双対であり、$\ft_c$は$T_c$のリー環で、$W$はワイル群である。前節にでてきた$T^*T^\vee/W$と同じものである。
Coulomb枝では、切断は$0$であり、$(\RR^3\otimes\ft_c)$は、ここで省略した場の成分から来るものである。$T^\vee_c$成分は接続から来るのだが、物理で`双対'とよんでいる、無限次元の接続の空間でのFourier変換を取るために、双対トーラス$T_c^\vee$に値をとり、かつスカラーになっている。この議論は、そのまま数学にのせるのは難しいと思うが、\subsecref{subsec:torus}と定理\ref{thm:classical}で見るように数学的に厳密な定義から出発して、$T^*T^\vee/W$ を再現することができ、なぜ双対トーラスになるのかも説明される。

$\cM_C$と$\cM_H$、より一般に真空のモジュライ空間は、超対称ゲージ理論の重要な情報を含んでおり、物理的にはゲージ理論を解析する上で、これを理解することは大切なステップである。特に、最初に与えた超対称性ゲージ理論が、真空のモジュライ空間をtargetにするような写像にいろいろな場を足して定められる超対称性場の理論と、低エネルギーにおいて等価になる。(ここで出てくる超対称性場の理論は、トポロジカル捻りをするとRozansky-Witten不変量を与えるものである。)

しかし、古典解に対応するような、ラグランジアンの最小値だけを見ていて、量子的な効果を含んだ場の理論の等価性を導くのは、過度な期待である。
物理学者は、そこでCoulomb枝は量子補正を受ける、と主張する。すなわち、Coulomb枝が$(T_c^\vee\times(\RR^3\otimes\ft_c))/W$であるのは古典的な記述であって、量子的な効果を受けたあとのCoulomb枝は、変更される、と主張する。ただし、超対称性から超ケーラー多様体であることは量子効果のあとも保たれる。
この変更が、超ケーラー構造の存在以外にどの程度$(T_c^\vee\times(\RR^3\otimes\ft_c))/W$を変更するのか、筆者には想像ができないが、数学的な定義のもとでは、$\cM_C$は$(T_c^\vee\times(\RR^3\otimes\ft_c))/W$と双有理同値であり、確かに変更していると取れなくもない。

というわけで、物理学者による$\cM_C$の定義は、$\cM_H$とは違って数学的には厳密とはいえず、そのままでは数学的に取り扱うことができない。筆者は、1996年11月にケンブリッジのニュートン研究所に滞在中に、Wittenの連続講演で初めてCoulomb枝の説明を聞いたが、研究対象として扱うことは長らくできなかった。出てくる超ケーラー多様体はよく知っているものであったので、頭の隅にずっと置いていたが、解決するのは難しいと考えていた。

新しい着想を得たのは、2014年秋にウォーリックでHananyの講演を聞いたときである。\linebreak[3]Hananyは、$\cM_C$の座標環$\CC[\cM_C]$の$\CC^\times$作用に関する指標を与える一般的な公式(モノポール公式)が成立すると説明した。この公式は、$G$のコウェイトに関する和で与えられ、足される項はコウェイトで定める具体的な式である。そして、知られているCoulomb枝の多くの例で、モノポール公式が確かに成立していることが、確かめられていた。

そこで、モノポール公式を再現するような空間を実現するためにはどうしたらいいかを逆に考えて発見したのが、\cite{2015arXiv150303676N}であり、その修正版の\cite{main}である。私が、どのように試行錯誤したかは\cite{2015arXiv150303676N}に説明したので、興味ある読者は参照されるとよいだろう。特に、$3$次元の位相的場の理論があると仮想して試行をしているところは、今後の発展の手がかりになるはずである、と期待しているところである。

\section{数学的な定義}\label{sec:def}

この節以降は、$G$は複素簡約群、$\bN$はその有限次元表現とする。$\bN$は既約でなくてもよく、$0$であってもよい。最初の節で述べた$\bM$は$\bN\oplus\bN^*$として与えられるが、ここから先は$\bM$は少なくとも表面上は出てこない。

$D = \Spec\CC[[z]]$をformal disk、$D^\times = \Spec\CC((z))$をformal punctured diskとする。$\bN((z))$, $\bN[[z]]$をそれぞれ $\bN_\cK$, $\bN_\cO$ で表す。同様に$G_\cK = G((z))$, $G_\cO = G[[z]]$とする。

アファイン・グラスマン$\Gr_G$は、モジュライ空間
\begin{equation*}
  \left.\left\{ (\scP,\varphi) \middle|
    \begin{aligned}[m]
      & \text{$\scP$は$D$上の(代数的な)$G$-主束}\\
      & \text{$\varphi\colon\scP|_{D^\times}\to G\times D^\times$は、$\scP$の$D^\times$上での自明化}
    \end{aligned}
   \right\}\middle/\text{isom.}\right.
\end{equation*}
として定義される。射影多様体の直極限としてのind-schemeの構造を持つことが知られている。集合論的には$\Gr_G = G_\cK/G_\cO$ と表される。すなわち、$\scP$の$D$上での自明化をとって、$\varphi$を$G_\cK$の元で表し、最初の自明化のambiguityの分の$G_\cO$で割って、$G_\cK/G_\cO$となる。

さらに、これに表現$\bN$に付随したベクトル束$\scP\times_G\bN$の切断$s$を付け加えた三つ組$(\scP,\varphi,s)$のモジュライ空間を$\cT$で表す。集合論的には$G_\cK\times_{G_\cO}\bN_\cO$である。$s$の展開を途中で止めることによって、$\cT$は射影多様体上のベクトル束の逆極限の直極限になる。
以下では、$\cT$や、その閉部分多様体のホモロジー群を取り扱うが、厳密には有限次元の空間のホモロジー群の極限として取り扱われる。

$\cT$の閉部分多様体$\cR$として、$\varphi(s)$が$D$まで伸びるという条件を課して、定められる空間と定義する。
\begin{equation*}
  \cR = \{ (\scP,\varphi,s) \mid \varphi(s)\in\bN_\cO \}/\text{isom.}
\end{equation*}
$\varphi$は$D^\times$上の自明化でしかないから、$\varphi(s)$は一般には原点に極を持つ有理型切断であって、その特異部分が$0$であるという条件を課したものが$\cR$である。集合論的には、
\(
   \cR = \{ [g,s]\in G_\cK\times_{G_\cO}\bN \mid gs\in\bN_\cO \}
\)
と記述できる。

この空間$\cR$が主要登場人物である。その意味は、より大きな空間
\begin{equation*}
  \{ (\scP_1, \varphi_1, s_1, \scP_2, \varphi_2, s_2)
  \in \cT\times\cT \mid \varphi_1(s_1) = \varphi_2(s_2) \}/\text{isom.}
\end{equation*}
を考えると、分かりやすいだろう。これは、$D$上の$G$主束と$D^\times$上の自明化および$\bN$に付随したベクトル束の切断の組が二つあって、切断が$D^\times$上で自明化を通じて等しい、というファイバー積$\cT\times_{\bN_\cK}\cT$に他ならない。$(\scP_2,\varphi_2)$が$\Gr_G$の原点、すなわち$\varphi_2$が$D$上の自明化に伸びているもの、になっているものが$\cR$に他ならない。
逆に、$\cR$への$G_\cO$の作用を用いて$G_\cK\times_{G_\cO}\cR$を考えると、これが上で出てきた空間$\cT\times_{\bN_\cK}\cT$に他ならない。

ゲージ理論的な視点では、$\cT\times_{\bN_\cK}\cT$は、二次元空間の上にある接続と切断の組が、原点のまわりでひねられている様子をあらわす空間である。$(\scP_1,\varphi_1,s_1)$がひねられる前で、$(\scP_2,\varphi_2,s_2)$がひねられる後であり、原点でひねられるだけなので、原点の外では一致している。本来は$3$次元のゲージ理論であるが、時間方向の動きは見ずに、前後の二つの瞬間だけを切り取って比べているので、$2$次元の記述になっている。

空間を準備し、次に$\cR$の$G_\cO$-同変Borel-Mooreホモロジー群$H^{G_\cO}_*(\cR)$を考える。厳密には、$\cT$の原点におけるファイバーの基本類が次数$0$になるように、次数をうまく定義する必要があるが、この点の詳細は略す。また、奇数次のホモロジーが消えていること、$H^*_G(\mathrm{pt})$上自由な加群になっていることなどは、アファイン・グラスマン多様体のSchubert胞体分割を考えると、ただちに従う。

$H^{G_\cO}_*(\cR)$には合成積
\begin{equation*}
  \ast\colon H^{G_\cO}_*(\cR)\otimes H^{G_\cO}_*(\cR)
  \to H^{G_\cO}_*(\cR)
\end{equation*}
が定義される。詳しい定義は、技術的なのでここでは略す。同変ホモロジー群のinduction
$H^{G_\cK}_*(\cT\times_{\bN_\cK}\cT)\cong H^{G_\cO}_*(\cR)$が有限次元の空間のときと同様に成り立っていると仮想的に考えて、さらに$\cT$は非特異であるとすると、通常の合成積の定義が、$(i,j)$成分への射影
\begin{equation*}
  \cT\times_{\bN_\cK}\cT\times_{\bN_\cK}\cT \xrightarrow{p_{ij}}
  \cT\times_{\bN_\cK}\cT \qquad (i,j) = (1,2), (2,3), (1,3)
\end{equation*}
を用いて
\begin{equation*}
  c\ast c' = p_{13*}(p_{12}^*c \cap p_{23}^*c')
\end{equation*}
と定義される。$H^{G_\cK}_*(\cT\times_{\bN_\cK}\cT)$ が定義されるかどうかは不明であり、$\cT$は非特異でないので、このままの定義がうまくできているのかどうかは分からないが、実際には$H^{G_\cO}_*(\cR)$の上に合成積$\ast$が定義される。

このとき次が成立する。
\begin{Theorem}\rm
  $(H^{G_\cO}_*(\cR),\ast)$ は可換環である。
\end{Theorem}

合成積で環を構成する手法は、幾何学的表現論で広く使われており、ワイル群の群環がSteinberg多様体から作られること、Kac-Moody Lie環の普遍展開環が箙多様体におけるSteinberg多様体の類似物から作られることなどが知られている。これらの例では得られるものは、非可換環であり、合成積の一般論からは$\ast$が可換になる理由はなく、上の定理は今の状況の特殊性を表している。

ただし、幾何学的佐武対応を思い起こせば、可換性は不思議ではない。幾何学的佐武対応では、アファイン・グラスマン$\Gr_G$上の$G_\cO$-同変な偏屈層のなすアーベル圏を考え、その上に合成積によってテンソル圏の構造を導入し、これが$G$のLanglands双対の有限次元表現の全体のなすテンソル圏と同値であることを主張する。後者のテンソル圏は可換、すなわち $V\otimes W\cong W\otimes V$ であるので、前者もそうである。この同型を幾何学的に説明するのがBeilinson-Drinfeldによるアファイン・グラスマンの1-パラメータ変形であり、これを使って上の定理が証明される。(論文では、計算による直接証明も与えている。)

さて、$(H^{G_\cO}_*(\cR),\ast)$は可換環になったので、そのスペクトラムとしてアファイン多様体を導入することができる。これが、Coulomb枝の数学的な定義である。
\begin{equation*}
   \cM_C = \Spec (H^{G_\cO}_*(\cR),\ast)
\end{equation*}
さらに、$(H^{G_\cO}_*(\cR),\ast)$が有限生成であることや、整であることを証明できるので、$\cM_C$は既約なアファイン多様体である。また正規であることも示されている。

変形量子化は、次のようにして与えられる。formal disk $D$に、$\CC^\times$がloop rotation$z\mapsto tz$により作用する。この作用は、今まで使ってきた様々な空間への作用を引き起こす。特に、$G_\cO$に作用して、半直積$G_\cO\rtimes\CC^\times$を考えることができ、$\cR$に$G_\cO\rtimes\CC^\times$が作用する。そこで、同変Borel-Mooreホモロジー群$H^{G_\cO\rtimes\CC^\times}(\cR)$を考え、合成積を同じように導入する。こうして{\bf 量子化されたCoulomb枝}を
\begin{equation*}
   \cAh = (H^{G_\cO\rtimes\CC^\times}(\cR),\ast)
\end{equation*}
と定義する。

なお、アファイン・グラスマン多様体や、その類似の合成積を考えるのは、以前から\cite{MR3013034,MR2135527,MR2422266}で考えられており、合成積の定義を厳密に書き下す際には、それを参考にした。\cite{MR3013034}では、アファイン・グラスマンの代わりにアファイン旗多様体、同変ホモロジー群の代わりに同変K群が用いられているが、$\bN = \mathfrak g$の場合を扱っていると思ってよい。出てくる代数は、Cherednikの二重アファイン・ヘッケ代数(DAHA)である。Coulomb枝のようにアファイン・グラスマン多様体にすれば、そのspherical partになり、同変ホモロジー群になれば楕円版の代わりに三角関数版のDAHAになる。対応するCoulomb枝は$\ft\times T^\vee/W$であり、量子補正がない、ということになる。

\cite{MR2135527,MR2422266}では、$\bN=0$の場合を取り扱っている。出てくるものは、$G$のLanglands双対の戸田格子であるが、詳細は略す。

\section{例}

前節の構成は、無限次元空間のホモロジーを使うもので、ずいぶんと回りくどい構成に見えるかもしれないので、簡単な例をあげよう。

\subsection{}\label{subsec:torus}

$G = \CC^\times$とし、$\bN = 0$ とする。これは、一番自明な例である。$\bN=0$なので、$\cR$はアファイン・グラスマン$\Gr_G$に他ならず、また$G=\CC^\times$なので、$\Gr_G$は$D$上の直線束とその$D^\times$上での自明化のモジュライ空間に他ならない。reduced schemeを取ると、$\Gr_G$は整数$\ZZ$でパラメトライズされた離散的な空間になる。実際、$\varphi(z) = z^n$ ($n\in\ZZ$)が対応する点をあらわす。よって
\begin{equation*}
   H^{G_\cO}_*(\cR) = \bigoplus_n H^{\CC^\times}_*(\mathrm{pt})
\end{equation*}
となる。$H^{\CC^\times}_*(\mathrm{pt})$は、一変数の多項式環$\CC[w]$である。これが各整数$n$の上に乗っているので、$m$の上の多項式と$n$の上の多項式を掛けるとどうなるかを、合成積の定義に戻って計算する。合成積の定義を説明しなかったので、チェックすることはできないが、$G=\CC^\times$の場合には、テンソル積を取る写像
\begin{equation*}
   \Gr_{\CC^\times}\times\Gr_{\CC^\times}\xrightarrow{\otimes}\Gr_{\CC^\times}
\end{equation*}
があり、これがホモロジー群に引き起こすpushforward凖同型が$\ast$に他ならない。すると$m$の上の$f(w)$と$n$の上の$g(w)$を掛けたものは、$m+n$の上の$f(w)g(w)$になる。すなわち、$n=1$ の上の$1$(基本類に対応する)を$x$であらわすと、
\begin{equation*}
   H^{G_\cO}_*(\cR) \cong \CC[w,x^\pm] = \CC[\CC\times\CC^\times]
\end{equation*}
となる。したがって、今のばあいのCoulomb枝は$\CC\times\CC^\times$である。これは$\RR^3\times S^1$であるから、この場合のCoulomb枝は量子補正を受けないことを意味しており、ゲージ理論が自明なことの反映である。

もう一歩、精密に見るために$G$はトーラス$T$で、表現はやはり$0$であるとする。$\Gr_T$は離散的な空間で、$\Hom(\CC^\times,T)$でパラメトライズされている。従って、$H^{T_\cO}_*(\cR) = \bigoplus_{\lambda\in\Hom(\CC^\times,T)} H^*_T(\mathrm{pt})$である。$H^*_T(\mathrm{pt})$は、$T$のLie環$\ft$上の多項式環$\CC[\ft]$である。一方、$\lambda$に対応する元を$e^\lambda$と書くと、上と同様に$e^\lambda\ast e^\mu = e^{\lambda+\mu}$ となる。これは、$T$の双対 $T^\vee$ の指標 ($\Hom(T^\vee,\CC^\times) = \Hom(\CC^\times, T)$) と見なすことができるから、Coulomb枝は$\ft\times T^\vee = T^* T^\vee$である。

\subsection{}\label{subsec:C2}

次に$G$は$\CC^\times$のままで, 表現を$\bN = \CC$ と標準表現に取ろう。
% このとき$\bN_\cO = \CC[z]$である。
$\Gr_{\CC^\times}$は上で説明したように$\ZZ$でパラメトライズされる離散的な空間であり、$\cR$は各整数$n$の上にベクトル空間が乗っているものである。条件は$\varphi(z) = z^n$によって原点に特異点が生じないというものであるから、
\begin{equation*}
  \cR = \bigsqcup_{n\in\ZZ} z^{n}\CC[z]\cap \CC[z] =
  \bigsqcup_{n\in\ZZ} z^{\max(0,n)}\CC[z]
\end{equation*}
である。各整数の上に乗っているものは、ベクトル空間でありThom同型により$H^{G_\cO}(\cR) \cong \bigoplus_n H^{\CC^\times}_*(\mathrm{pt})$となる。すなわちベクトル空間としては、上の例と同じである。しかし合成積は、上の例とは$n>0$の上のホモロジー類と$n<0$の上のホモロジー類の積が変わってくる。定義を省略したので、最後のポイントだけいうと、$n=1$の基本類と$n=-1$の基本類を掛けたものが
\begin{equation*}
  z\CC[z] \to \CC[z]
\end{equation*}
の押し出し写像による、基本類の像になる。これは余次元$1$の部分空間であるから、同変ホモロジー群の元としては、$w$を基本類に掛けたものになる。したがって$n=1$の基本類を$x$, $n=-1$の基本類を$y$とすると、$xy = w$が成り立つ。この計算から
\begin{equation*}
  H^{G_\cO}_*(\cR) \cong \CC[w,x,y]/(w=xy) \cong \CC[x,y] = \CC[\CC^2]
\end{equation*}
が従う。よって今の場合のCoulomb枝は$\CC^2$である。

表現をウェイトが$N$の一次元表現に取り替えると、最後の部分の計算が$z^{|N|}\CC[z]\to \CC[z]$のpushforwardに置き換わり、座標環は$\CC[w,x,y]/(w^{|N|}=xy)$となる。これは、$A_{|N|-1}$型の単純特異点に他ならない。

\section{いくつかの構造}

この節では、Coulomb枝$\cM_C$が持つ構造について解説する。いずれも物理的には発見されていたが、これが\secref{sec:def}の定義のもとで数学的に厳密に実現されることがポイントである。

\subsection{}

$H^{G_\cO}_*(\cR)$は、ホモロジーの次数の半分により次数付けられた環になる。つまり、$\CC[\cM_C] = \bigoplus_d \CC[\cM_C]_d$ と分解し、$\CC[\cM_C]_d\cdot\CC[\cM_C]_{d'}\subset\CC[\cM_C]_{d+d'}$となる。これは、$\cM_C$に$\CC^\times$の作用が与えられていることを意味する。実際、$\CC[\cM_C]_d$は、$\CC^\times$がウェイト$d$で作用するウェイト空間である。

上の例では、$\CC\times\CC^\times$と$\CC^2 = \CC\times\CC$のそれぞれ第一成分への標準的な作用になっている。(正確には、後者は$x$がウェイト$1$で、$y$がウェイト$0$である。)

なお、次数の定義を省いたので説明が不足しているが、次数は非負とは限らず、一般にはすべての整数の値を取りうる。従って、$\cM_C$は一般的には錘であるとは限らない。ここで、$\cM_C$が錘であるとは、$\CC[\cM_C]_d = 0$ ($d < 0$), $\CC[\cM_C]_0 = \CC$が成り立つときをいう。

物理的には、この$\CC^\times$作用は、$\SU(2)$作用の$S^1$への制限から来るものに、ある修正のあと一致すると期待されている。修正については説明しないが、次に述べる$\cM_C$へのハミルトニアンなトーラス作用を適当に組み合わせるものである。特に、$G$が半単純のときには、修正は必要ない。$\SU(2)$作用は、超ケーラー構造$I$, $J$, $K$がなす二次元球面 $S^2 = \{ a I + bJ + cK \mid a^2 + b^2 + c^2 = 1\}$ に$\SU(2)\to\SO(3)$を通じて標準的に作用するもので、一つの複素構造$I$を固定すると、$I$を保つ$S^1$の作用しか見えない。現在のところ、超ケーラー構造の定義は与えられていないので、$\SU(2)$の作用を数学的に与えることはできていないが、制限の$S^1$だけが見えている。

上の例では、$\CC\times\CC^\times = \RR^3\times S^1$として$\RR^3$を$\su(2)$と見れば、確かに$\SU(2)$の作用がある。ウェイトは半分になっている。$\CC^2$の場合は、$x$, $y$がそれぞれウェイト$-1/2$, $1/2$のハミルトニアンな$S^1$作用と合わせて、ウェイトが共に$1/2$の作用に直せば、やはりウェイトが半分になっていることを除き、$\CC^2$を四元数体$\mathbb H$とみた$\SU(2) = \grpSp(1)$の作用に一致する。(複素線形でないので、$\SU(2)$の$\CC^2$への標準的な表現とは異なり、四元数の右掛け算と左掛け算の違いがある。)

\subsection{}

$H^{G_\cO}_*(\cR)$は同変ホモロジー群であるから、$H^*_{G_\cO}(\mathrm{pt})\cong H^*_G(\mathrm{pt})$からの準同型を持つ。(ただし、$H^*_G(\mathrm{pt})$上の代数ではなく、合成積 $c\ast c'$は第二成分$c'$に関しては、自然には$H^*_G(\mathrm{pt})$-線形にはならず、変形量子化したものについては、確かに線形でない。)

これのスペクトラムを取ると、
\begin{equation*}
    \varpi\colon\cM_C\to \Spec H^*_G(\mathrm{pt})
\end{equation*}
を得る。よく知られているように、
\begin{equation*}
     H^*_G(\mathrm{pt}) = \CC[\operatorname{Lie}G]^G
    = \CC[\ft]^W
\end{equation*}
であるから、$\Spec H^*_G(\mathrm{pt}) = \ft/W$であり、これはアファイン空間である。ここで、$\ft = \operatorname{Lie}T$である。

この構成は、変形量子化したあとも残り、
\begin{equation*}
    H^{G\times\CC^\times}_*(\mathrm{pt})\to
    \cAh = H^{G_\cO\rtimes\CC^\times}(\cR)
\end{equation*}
という単射な環準同型がある。これは、変形量子化が大きな可換環を含んでいることを意味しており、これの $\hbar = 0$を考えることにより、$\varpi$はポアソン可換であることを導く。すなわち、$\operatorname{Lie}T/W$上の関数 $f$, $g$ を $\varpi$で引き戻したものは、Poisson可換である：$\{ \varpi^*f, \varpi^*g\} = 0$.

さらに次が成立する。
\begin{Theorem}\label{thm:classical}\rm
$\varpi$のgenericなファイバーは、$T^\vee$である。より強く、次の可換図式がある。上の横矢印は双有理写像である。
\begin{equation*}
    \xymatrix{
    \cM_C \ar@{.>}[rr] \ar[dr]_\varpi &&
    T^* T^\vee/W = \ft\times T^\vee/W \ar[dl]^{\text{第一射影}}
    \\
    & \ft/W &}
\end{equation*}
\end{Theorem}

これは、同変ホモロジー群の局所化定理の帰結である。局所化定理は、$H^*_T(\mathrm{pt})$の商体を$\mathbb F$とするとき、
\begin{equation*}
    H^{T_\cO}_*(\cR)\otimes_{H^*_T(\mathrm{pt})}\mathbb F
    \cong H^{T_\cO}_*(\cR^T)\otimes_{H^*_T(\mathrm{pt})}\mathbb F
\end{equation*}
が成り立つという主張である。ここで、$\cR^T$は$\cR$の$T$-固定点の集合であり、同型写像は、包含写像$\cR^T\hookrightarrow\cR$のpushforward準同型である。これと、$H^{G_\cO}_*(\cR)$は、$H^{T_\cO}_*(\cR)$の$W$-不変部分であるという事実を組み合わせると、$\cR^T$の同変ホモロジー群を決定すれば良いことになるが、$\cR^T$が$\Gr_T\times\bN^T$であることと、\subsecref{subsec:torus}の計算から、$\ft\times T^\vee$であることが分かる。

操作$\otimes_{H^*_T(\mathrm{pt})}\mathbb F$は、$\ft/W$のgeneric pointに制限することであり、同変ホモロジー群を考えることが$\ft/W$上の族を考えるという幾何学的な描像に対応しているという、よく知られた哲学の有効性をあらわす典型的な議論である。

以上で、$\varpi$はポアソン交換していて、ファイバーが代数的トーラスであることから、$\varpi$はLiovilleの意味で可積分系である。変形量子化$\cAh$はその量子化である。

\subsection{}

アファイン・グラスマン $\Gr_G$は位相的には基点付きループ群 $\Omega G$であることが知られており、特にその連結成分は$G$の基本群 $\pi_1(G)$に一致する。ホモロジー群は連結成分に応じて分解するが、これは合成積とcompatibleである。すなわち$\gamma\in\pi_1(G)$に対応する$\cR$の連結成分を$\cR_\gamma$と書くとき、$H^{G_\cO}_*(\cR_\gamma)\ast H^{G_\cO}_*(\cR_{\gamma'})\subset
H^{G_\cO}_*(\cR_{\gamma+\gamma'})$となる。($\pi_1(G)$が可換であることはよく知られている。) 従って $H^{G_\cO}_*(\cR)$は、$\pi_1(G)$で次数付けられた環である。

これを $\cM_C = \Spec H^{G_\cO}_*(\cR)$側で考えると、$\pi_1(G)$のポントリャーギン双対$\Hom(\pi_1(G),\linebreak[3]\CC^\times)$が$\cM_C$に作用することになる。たとえば、上の例では$\pi_1(G) = \pi_1(\CC^\times) = \ZZ$であり、ポントリャーギン双対は$\CC^\times$である。$\CC^\times$作用は、最初の例では$\CC\times\CC^\times$の第二成分への自然な作用であり、二番目の例の場合は$x$がウェイト$1$で、$y$がウェイト$-1$である。

この作用は、変形量子化 $H^{G_\cO\rtimes\CC^\times}_*(\cR)$にも自然に伸びていることから、シンプレクティック形式を保っていることも従う。

$G$が半単純のときには、$\pi_1(G)$は有限群で、そのポントリャーギン双対も有限群になってしまうが、トーラスが現れるのは$\Hom(G,\CC^\times)$が自明でない場合である。このとき、運動量写像は$\varpi$に$\operatorname{Lie}G\to\operatorname{Lie}\CC^\times$を合成したもので与えられ、特に作用はハミルトニアンである。

\subsection{}

$\bN$が、$G$を正規部分群として含む大きな群$\widetilde G$の表現の制限として現れている場合を考える。物理では、商群 $\widetilde G/G$はフレーバー対称性の群とよばれる。これを$G_F$で表わす。

$\widetilde G_\cO$は$\cR$に作用するので、大きな群の同変ホモロジー$H^{\widetilde G_\cO}_*(\cR)$を考えることができる。合成積により、$H^*_{G_F}(\mathrm{pt})$上の代数になり、対応するスペクトラムは、$\Spec H^*_{G_F}(\mathrm{pt}) = \CC[\operatorname{Lie}G_F]^{G_F}$上の多様体の族になり、原点のファイバーが元の$\cM_C$である。すなわち、$\cM_C$は$\operatorname{Lie}G_F\dslash G_F$でパラメトライズされた変形を持つ。

また、説明は省略するが、変形に対応するような$\cM_C$の部分特異点解消(の候補)を構成することもできる。

この節と前節では、$\Hom(G,\CC^\times)$と$G_F$が$\cM_C$に導く構造を調べたが、Higgs枝$\cM_H$に導く性質を考えることは有用である。まず、$G_F$であるが、$\cM_H = \bM\tslash G$であるから、$G_F$は$\cM_F$に作用する。
一方で、$\Hom(G,\CC^\times)$があると、対応する$\zeta\in\Hom(\operatorname{Lie}G,\operatorname{Lie}\CC^\times)$を考えて、運動量写像のレベル集合を$\mu=0$から$\mu=\zeta$に変形することができる。
すなわち、$\Hom(G,\CC^\times)$と$G_F$が$\cM_C$と\linebreak[3]$\cM_H$に誘導する構造は、それぞれ群作用と変形であるが、両者は$\cM_C$と$\cM_H$で入れ替わっている。

\subsection{}

前節と、前々節の構造の例として、トーリック超ケーラー多様体を考える。こ
れには、トーラスの完全列
\begin{equation*}
    1 \to T = (\CC^\times)^{d-n} \to \widetilde T = (\CC^\times)^d
    \to T_F = (\CC^\times)^n \to 1
\end{equation*}
が与えられたとする。$\widetilde T$の標準的な表現$\bN = \CC^d$を取り、その$T$への制限も$\bN$で表わす。さて、$\cM_C(\widetilde T,\bN)$は
\subsecref{subsec:C2}の計算より$\CC^{2d}$となる。前々節の構成により$\pi_1(\widetilde T)$のポントリャーギン双対が$\CC^{2d}$に作用するが、これは$\widetilde T$の双対トーラス$\widetilde T^\vee$に他ならない。$T_F$の双対トーラスは$T_F^\vee$はその部分トーラスであり、前々節の構成をもう少し進めると $T$に関するCoulomb枝 $\cM_C(T,\bN)$は、$\CC^{2d}$の$T_F^\vee$に関するシンプレクティック商
$\CC^{2d}\tslash T_F^\vee$に他ならない。これは、双対トーラスの完全列
\begin{equation*}
    1\to T_F^\vee \to \widetilde T^\vee \to T^\vee \to 1
\end{equation*}
を考えて、$\widetilde T^\vee$の表現$\bM = \CC^d\oplus(\CC^d)^*$に関するHiggs枝ということもできる。すなわち、$T$と$T_F^\vee$を入れ替えることによって、Higgs枝とCoulomb枝が入れ替わっている。

\section{箙ゲージ理論}

現在のところ、Higgs枝が箙多様体になるような$(G,\bN)$に対応するCoulomb枝についてが一番よく調べられている。$Q$ を箙とし、$Q_0$をその頂点の集合、$Q_1$を向きの付けられた辺の集合とする。$h\in Q_1$に対し、その始点と終点を$\vout{h}$, $\vin{h}$で表わす。二つの$Q_0$で次数付けられたベクトル空間 $V = \bigoplus V_i$, $W = \bigoplus W_i$が与えられたとき、
\begin{equation*}
    \begin{split}
        & G = \prod_{i\in Q_0}\GL(V_i),
        \\
        & \bN = \bigoplus_{h\in Q_1}\Hom(V_{\vout{h}},V_{\vin{h}})
        \oplus\bigoplus_{i\in Q_0} \Hom(W_i, V_i)
    \end{split}
\end{equation*}
が、箙ゲージ理論である。ただし、$G$の$\bN$への作用は自然なものである。

$Q$が$ADE$型の場合には、$\cM_C$は原点に特異点を持った$\RR^3$の上のモノポールのモジュライ空間になると物理的には洞察されていたが、この空間の代数幾何的な対応物が、先の数学的な定義のもとで示されている。(\cite{2016arXiv160403625B}) ここで、モノポールの構造群は、$Q$に対応する(adjoint型の)複素単純リー群$G_Q$であり、$V$の次元は、モノポールの次数に対応し、$W$の次元は特異点の情報を与える。
代数幾何的な対応物は、一般の場合は記述は簡単ではないが、$\mu = \sum \dim W_i \varpi_i - \dim V_i \alpha_i$が支配的なときには、$G_Q$のアファイン・グラスマンを考え、$\lambda = \sum \dim W_i \varpi_i$と$\mu$に対応するSchubert多様体$\overline{\Gr}_{G_Q}^{\la}$, $\overline{\Gr}_{G_Q}^{\mu}$を取って、$\overline{\Gr}_{G_Q}^\mu$の横断切片と$\overline{\Gr}_{G_Q}^{\la}$の交わりが$\cM_C$である。

幾何学的佐竹対応によってアファイン・グラスマンは $G_Q$のラングランズ双対の表現論と結びついていたが、一方で箙多様体のホモロジー群には、$G_Q$の表現の構造が入っていた。始めに述べたシンプレクティック双対性は、この二つの構成が`双対'によって結びついていることを主張するように定式化される(べきである)。

この結果の証明のためには、$\cM_C$を決定する次のような処方箋を用いる。
\begin{enumerate}
      \item 
    まず、$\cM_C$の候補になる空間を作る。これは、多くの場合は、物理学者の答えを採用する。
      \item
    次に、その候補の空間に、$\varpi$に対応すると期待される可積分系を作る。

      \item
    その可積分系が平坦な族であること、$\cM_C$が正規であることをチェックする。

      \item 
    
$\cM_C$とその候補の間の $T^* T^\vee/W$を通じた双有理写像が$\ft/W$の余次元$2$の集合を除いて拡張することをチェックする。
\end{enumerate}
最後の余次元$2$の集合を除けば十分であるところは、正規性の帰結である。
同変ホモロジーの局所化定理の応用で、genericには$T^*T^\vee$であることを説明したが、余次元$1$のところも同様の議論で、階数$1$の群のCoulomb枝を決定する問題に帰着できる。階数$1$の場合は、Coulomb枝は$\CC^3$の超曲面として実現できることが示されており、決定されている。したがって、(4)は、易しいステップである。現状では、(3)を示す部分が、ケースバイケースで行われていてキーポイントになっている。

アファイン$ADE$型の場合は、有限型のモノポールの代わりにインスタントンを考えればよい。ただし、$\RR^4$上のインスタントンではなく、Taub-NUT空間上のインスタントンにするのが正確なので、微妙な問題があり、特に上でいうところの$\mu$が支配的な場合は$\RR^4$上でもTaub-NUT空間上でも、複素シンプレクティック多様体としては変わらないと期待されている。

インスタントンのモジュライ空間については、(3)の性質が証明されていないので、現在のところCoulomb枝の決定までは至っていない。

(3)は微妙な性質である。例えばべき零軌道は、$A$型のときは常に正規であるが、一般にはそうでない。一方、Coulomb枝は常に正規である。古典型のべき零軌道やそのSlodowy横断切片は、Higgs枝として現れることが知られているので、対応するCoulomb枝も、そうなっていると安直には考えられるが、正規性の問題から、そうそう単純ではなさそうである。Hananyらは、正規化を取ればよいと考えているようであるが、まだまだ十分な根拠があるとはいえないのではないだろうか？

アファイン$A$型のときには、インスタントンのモジュライ空間を直接取り扱う代わりに、Cherkisの弓箭多様体(bow varietyの和訳)を用いる。弓箭多様体は、Nahm方程式とよばれる非線形常微分方程式の解を用いてあらわされているので、そのままでは取扱いにくいが、\cite{2016arXiv160602002N}により、箙多様体の変種として書き直し、(3)の性質を証明した。したがって、アファイン$A$型の箙ゲージ理論のCoulomb枝は決定された。

\section{量子化されたCoulomb枝}

量子化されたCoulomb枝$\cAh$については、多様体の決定に比べると分かっている例は少ない。

前に、随伴表現$\bN=\mathfrak g$のときにspherical DAHAが出てくることを言及したが、$G=\GL(k)$のときは、ジョルダン箙に対応する箙ゲージ理論で、$V=\CC^k$, $W=0$の場合であると思うことができる。これを一般化して$V=\CC^k$, $W=\CC^r$と変えると、$\cAh$はwreath積$\ZZ/r\ZZ\wr S_k = (\ZZ/r\ZZ)^k\rtimes S_k$の有理Cherednik代数のspherical partになる。\cite{2016arXiv160800875K} 対応するCoulomb枝は$\operatorname{Sym}^k\CC^2/(\ZZ/r\ZZ)$である。

有限$ADE$型の箙ゲージ理論の場合は、\cite{2016arXiv160403625B}のAppendixにおいて$\cAh$がshited Yangianとして同型であることが示された。ただし、前節で言及した$\mu$が支配的という条件を仮定した下で証明されており、一般の場合は未解決である。

\bibliographystyle{myamsalpha}
\bibliography{nakajima,mybib,coulomb}    

\def\cprime{$'$} \def\cprime{$'$} \def\cprime{$'$} \def\cprime{$'$}
  \def\cprime{$'$}
  \providecommand{\noopsort}[1]{}\def\cftil#1{\ifmmode\setbox7\hbox{$\accent"5E#1$}\else
  \setbox7\hbox{\accent"5E#1}\penalty 10000\relax\fi\raise 1\ht7
  \hbox{\lower1.15ex\hbox to 1\wd7{\hss\accent"7E\hss}}\penalty 10000
  \hskip-1\wd7\penalty 10000\box7}
\providecommand{\bysame}{\leavevmode\hbox to3em{\hrulefill}\thinspace}
\providecommand{\MR}{\relax\ifhmode\unskip\space\fi MR }
\providecommand{\MRhref}[2]{%
  \href{http://www.ams.org/mathscinet-getitem?mr=#1}{#2}
}
\providecommand{\href}[2]{#2}
\begin{thebibliography}{{\noopsort{01}}BFN16}

\bibitem[BF08]{MR2422266}
R.~Bezrukavnikov and M.~Finkelberg, \emph{Equivariant {S}atake category and
  {K}ostant-{W}hittaker reduction}, Mosc. Math. J. \textbf{8} (2008), no.~1,
  39--72, 183. \MR{2422266 (2009d:19008)}

\bibitem[BFM05]{MR2135527}
R.~Bezrukavnikov, M.~Finkelberg, and I.~Mirkovi{\'c}, \emph{Equivariant
  homology and {$K$}-theory of affine {G}rassmannians and {T}oda lattices},
  Compos. Math. \textbf{141} (2005), no.~3, 746--768. \MR{2135527
  (2006e:19005)}

\bibitem[BFN16a]{main}
A.~{\noopsort{01}}~Braverman, M.~Finkelberg, and H.~Nakajima, \emph{{Towards a
  mathematical definition of Coulomb branches of $3$-dimensional $\mathcal N=4$
  gauge theories, II}}, ArXiv e-prints (2016),
  \href{http://arxiv.org/abs/1601.03586}{{\ttfamily arXiv:1601.03586
  [math.RT]}}.

\bibitem[BFN16b]{2016arXiv160403625B}
\bysame, \emph{{Coulomb branches of
  $3d$ $\mathcal N=4$ quiver gauge theories and slices in the affine
  Grassmannian (with appendices by Alexander Braverman, Michael Finkelberg,
  Joel Kamnitzer, Ryosuke Kodera, Hiraku Nakajima, Ben Webster, and Alex
  Weekes)}}, ArXiv e-prints (2016),
  \href{http://arxiv.org/abs/1604.03625}{{\ttfamily arXiv:1604.03625
  [math.RT]}}.

\bibitem[BLPW14]{2014arXiv1407.0964B}
T.~{Braden}, A.~{Licata}, N.~{Proudfoot}, and B.~{Webster},
  \emph{{Quantizations of conical symplectic resolutions II: category $\mathcal
  O$ and symplectic duality}}, ArXiv e-prints (2014),
  \href{http://arxiv.org/abs/1407.0964}{{\ttfamily arXiv:1407.0964 [math.RT]}}.

\bibitem[KN16]{2016arXiv160800875K}
R.~{Kodera} and H.~{Nakajima}, \emph{{Quantized Coulomb branches of Jordan
  quiver gauge theories and cyclotomic rational Cherednik algebras}}, ArXiv
  e-prints (2016), \href{http://arxiv.org/abs/1608.00875}{{\ttfamily
  arXiv:1608.00875 [math.RT]}}.

\bibitem[Nak92]{MR1193019}
中島　啓, Einstein計量の収束定理とALE空間,
  数学 \textbf{44} (1992), no.~2, 133--146. \MR{MR1193019 (93k:53044)}

\bibitem[Nak16]{2015arXiv150303676N}
H.~{Nakajima}, \emph{{Towards a mathematical definition of Coulomb branches of
  $3$-dimensional $\mathcal N=4$ gauge theories, I}}, Adv. Theor. Math. Phys.
  \textbf{20} (2016), no.~3, 595--669,
  \href{http://arxiv.org/abs/1503.03676}{{\ttfamily arXiv:1503.03676
  [math-ph]}}.

\bibitem[NT16]{2016arXiv160602002N}
H.~{Nakajima} and Y.~{Takayama}, \emph{{Cherkis bow varieties and Coulomb
  branches of quiver gauge theories of affine type $A$}}, ArXiv e-prints
  (2016), \href{http://arxiv.org/abs/1606.02002}{{\ttfamily arXiv:1606.02002
  [math.RT]}}.

\bibitem[VV10]{MR3013034}
M.~Varagnolo and E.~Vasserot, \emph{Double affine {H}ecke algebras and affine
  flag manifolds, {I}}, Affine flag manifolds and principal bundles, Trends
  Math., Birkh\"auser/Springer Basel AG, Basel, 2010, pp.~233--289.
  \MR{3013034}

\end{thebibliography}
\end{CJK}

\end{document}